\documentclass[twoside,a4paper,12pt,centertags,reqno]{amsart} % 13pt displays smaller and better than 12pt
\usepackage{amsmath,amssymb,verbatim,vmargin}
\usepackage{color}
\usepackage{tikz}
\usepackage{tikz-cd}
\usetikzlibrary{matrix,calc}
\usepackage{hyperref}
\hypersetup{colorlinks,bookmarks=true,linktocpage=true,citecolor=blue,linkcolor=magenta}

\usepackage{eulervm}   % formulas round   

\allowdisplaybreaks % align break
\usepackage{color}
\usepackage[all]{xy} % diagrams

\usepackage{mathtools}
\usepackage{MnSymbol} %wedge righthalfcup

\theoremstyle{plain}
\newtheorem{thm}{Theorem}[section]
\newtheorem{lem}[thm]{Lemma}

\newtheorem{prop}[thm]{Proposition}

\theoremstyle{definition}
\newtheorem{defn}[thm]{Definition}
\newtheorem{rem}[thm]{Remark}

\usepackage{enumerate}
\usepackage{enumitem}

\newcommand{\bC}{{\mathbb C}}

\newcommand{\bR}{{\mathbb R}}

\newcommand{\wh}{\widehat}

\def\barint_#1{\mathchoice
            {\mathop{\vrule width 6pt
height 3 pt depth -2.5pt
                    \kern -9.5pt
\intop \kern -4pt}\nolimits_{#1}}%
            {\mathop{\vrule width 5pt height
3 pt depth -2.6pt
                    \kern -6.5pt
\intop \kern -4pt}\nolimits_{#1}}%
            {\mathop{\vrule width 5pt height
3 pt depth -2.6pt
                    \kern -6pt
\intop \kern -4pt}\nolimits_{#1}}%
            {\mathop{\vrule width 5pt height
3 pt depth -2.6pt
          \kern -6pt \intop \kern -4pt}\nolimits_{#1}}}
          
           \def\bariint_#1{\mathchoice
            {\mathop{\vrule width 15pt
height 3 pt depth -2.5pt
                    \kern -15.8pt
\intop \kern -8pt\intop \kern -4pt}\nolimits_{#1}}%
            {\mathop{\vrule width 9pt height
3 pt depth -2.6pt
                    \kern -10.5pt
\intop \kern -8pt\intop \kern -4pt}\nolimits_{#1}}%
            {\mathop{\vrule width 9pt height
3 pt depth -2.6pt
                    \kern -10pt
\intop \kern -8pt\intop \kern -4pt}\nolimits_{#1}}%
            {\mathop{\vrule width 9pt height
3 pt depth -2.6pt
          \kern -8pt \intop \kern -10pt\intop \kern -4pt}
      \nolimits_{  #1}}}

\def\barintlim_#1{\mathchoice
            {\mathop{\vrule width 6pt
height 3 pt depth -2.5pt
                    \kern -8.8pt
\intop \kern -4pt}\limits_{#1}}%
            {\mathop{\vrule width 5pt height
3 pt depth -2.6pt
                    \kern -6.5pt
\intop \kern -4pt}\limits_{#1}}%
            {\mathop{\vrule width 5pt height
3 pt depth -2.6pt
                    \kern -6pt
\intop \kern -4pt}\limits_{#1}}%
            {\mathop{\vrule width 5pt height
3 pt depth -2.6pt
          \kern -6pt \intop \kern -4pt}\limits_{#1}}}
          
           \def\bariintlim_#1{\mathchoice
            {\mathop{\vrule width 15pt
height 3 pt depth -2.5pt
                    \kern -15.8pt
\intop \kern -8pt\intop \kern -4pt}\limits_{#1}}%
            {\mathop{\vrule width 9pt height
3 pt depth -2.6pt
                    \kern -10.5pt
\intop \kern -8pt\intop \kern -4pt}\limits_{#1}}%
            {\mathop{\vrule width 9pt height
3 pt depth -2.6pt
                    \kern -10pt
\intop \kern -8pt\intop \kern -4pt}\limits_{#1}}%
            {\mathop{\vrule width 9pt height
3 pt depth -2.6pt
          \kern -8pt \intop \kern -10pt\intop \kern -4pt}
      \limits_{  #1}}}
          
\renewcommand{\iint}{\int \kern -8pt\int}       

\numberwithin{equation}{section}
\makeatletter
\@namedef{subjclassname@2020}{\textup{2020} Mathematics Subject Classification}
\makeatother

\title[On Ozawa-Rogers for Klein-Gordon]
{On the bilinear estimate of Ozawa and Rogers for the one-dimensional Klein-Gordon equation
and some related lower Jacobian estimates}

\author{Shirong Chen}
\address{School of Mathematics and Computer Sciences, Gannan Normal University, Ganzhou 341000, People's Republic of China}
\email{somxiaorong@163.com}

\author{Yi C. Huang}
\address{School of Mathematical Sciences, Nanjing Normal University, Nanjing 210023, People's Republic of China}
\email{Yi.Huang.Analysis@gmail.com}
\urladdr{https://orcid.org/0000-0002-1297-7674}

\author{Shaozhen Xu}
\address{School of Information Engineering, Nanjing Xiaozhuang University, Nanjing 211171, People's Republic of China}
\email{xushaozhen@njxzc.edu.cn}
 
\date{\today}


\begin{document}

\begin{abstract}
We give a natural convexity proof of an elementary inequality used by Ozawa and Rogers in proving their bilinear estimate for the one-dimensional Klein-Gordon equation.
This robust approach also enables us to derive the optimality of Ozawa-Rogers estimate and establish a new bilinear estimate. 
Our estimate is in sharp analogy with the bilinear estimate of Bez and Rogers for the wave equation.
We also present an alternative proof of some lower Jacobian estimates used crucially by Ozawa and Rogers in proving their Fourier restriction results on the whole hyperbola.
Based on interpolation, our proof is purely non-trigonometric.
\end{abstract}

\maketitle
    
\section{Introduction}

Denote by $\widehat{f}$ the Fourier transform of $f$, that is,
\begin{equation*}
    \widehat{f}(\xi)=\int_{\mathbb{R}} f(x)e^{-ix\xi} dx.
\end{equation*}
About ten years ago, Bez and Rogers obtained in \cite{BezRog13} the following sharp bilinear estimates for the wave equation on space-time $\mathbb{R}^{d+1}$ for $d\geq2$:
\begin{equation}\label{E1}
\begin{aligned}
    & \left\|e^{it\sqrt{-\Delta}}f_1e^{it\sqrt{-\Delta}}f_2 \right\|_{L_{t,x}^2(\mathbb{R}^{d+1})}^2\\ 
    &\qquad\leq C_d\int_{\mathbb{R}^{2d}} |\widehat{f_1}(\xi_1)|^2|\widehat{f_2}(\xi_2)|^2
    |\xi_1|^{\frac{d-1}{2}}|\xi_2|^{\frac{d-1}{2}} 
    \left(1-\frac{\xi_1 \cdot \xi_2}{|\xi_1||\xi_2|}\right)^{\frac{d-3}{2}} 
    d\xi_1 d\xi_2. 
\end{aligned}
\end{equation}
However, this estimate is no longer valid for $d=1$. 
For the Klein-Gordon equation, Ozawa and Rogers obtained in \cite{OzaRog14}\footnotemark
\footnotetext{Ozawa reported in El Escorial 2012 the main results obtained in this paper.} the following sharp bilinear estimate.

\begin{thm}[Ozawa-Rogers] \label{thmA}
Suppose that $\text{\rm supp}\,(\widehat{f_1}) \cap \text{\rm supp}\,(\widehat{f_2})=\emptyset$. Then
    \begin{align*}
         & \left\|e^{it\sqrt{1-\Delta}}f_1e^{it\sqrt{1-\Delta}}f_2 \right\|_{L_{t,x}^2(\mathbb{R}^{1+1})}^2 \\
        &\qquad\leq \frac{1}{(2\pi)^2} \int_{\mathbb{R}^{2}} |\widehat{f_1}(\xi_1)|^2|\widehat{f_2}(\xi_2)|^2 (1+\xi_1^2)^{3/4} (1+\xi_2^2)^{3/4} \frac{d\xi_1 d\xi_2}{|\xi_2-\xi_1|}.
    \end{align*}
\end{thm}

For related bilinear estimates in various settings, see Klainerman-Machedon \cite{KlaMac93, KlaMac96, KlaMac97},
Ozawa-Tsutsumi \cite{OzaTsu98}, Carneiro \cite{Car09}, Jeavons \cite{Jea14}, Bez-Jeavons-Ozawa \cite{BezJeaOza16}, Beltran-Vega \cite{BelVeg20} and the references therein.

The proof of Theorem \ref{thmA} relies crucially on the following elementary inequality.

\begin{lem}[Ozawa-Rogers] \label{lemA}
For $\xi_2 \geq \xi_1$,
\begin{equation}\label{E2}
    \frac{\xi_2-\xi_1}{(1+\xi_2^2)^{3/4}(1+\xi_1^2)^{3/4}} \leq 
    \frac{\xi_2}{\sqrt{1+\xi_2^2}} - \frac{\xi_1}{\sqrt{1+\xi_1^2}}.
\end{equation}
\end{lem}

This inequality was proved via ``an unlikely combination of five\footnotemark
\footnotetext{The Pythagorean, Difference, Double-angle, Product-to-sum, and Power-reduction identities.} trigonometric identities".
Here we propose a ``pedagogical" proof and illustrate its sharpness.
 
\begin{prop}\label{remark1}
The exponent $3/4$ in \eqref{E2} is optimal.
\end{prop}

Our proof of \eqref{E2} uses convexity in nature and involves straightforward computations of derivatives. 
It also motivates us to derive the following bilinear estimate.
 \begin{thm}\label{thm1}
     Suppose that $\text{\rm supp}\,(\widehat{f_1}) \cap \text{\rm supp}\,(\widehat{f_2})=\emptyset$. Then
    \begin{equation}\label{E3}
 \begin{aligned}
         & \left\|e^{it\sqrt{1-\Delta}}f_1e^{it\sqrt{1-\Delta}}f_2 \right\|_{L_{t,x}^2(\mathbb{R}^{1+1})}^2 \\
         &\qquad\leq \frac{1}{(2\pi)^2} \int_{\mathbb{R}^{2}} |\widehat{f_1}(\xi_1)|^2|\widehat{f_2}(\xi_2)|^2 \left(1-\frac{(1,\xi_1)\cdot(1,\xi_2)}{\sqrt{1+\xi_1^2}\sqrt{1+\xi_2^2}}\right)^{-1} d\xi_1 d\xi_2. 
         \end{aligned}
\end{equation}
 \end{thm}
 
Note that the weight factor in \eqref{E3} is in ``natural" analogy with \eqref{E1} for $d=1$ in the following sense: 
for the Klein-Gordon equation, the ``right" spherical distance to be talked about 
concerns the two unit vectors ${(1,\xi_1)}/{\sqrt{1+\xi_1^2}}$ and ${(1,\xi_2)}/{\sqrt{1+\xi_2^2}}$.

\begin{rem}
Related inequalities similar to \eqref{E2} were also obtained in \cite{OzaRog14} and crucially used for their Fourier restriction results. 
We discuss them in Section \ref{sec:restr}.
\end{rem}

\section{A new proof of Lemma \ref{lemA}}
Fix $\xi_1$ and consider the function 
\begin{equation*}
    f(x)=\left(\xi_1 ^2+1\right)^{3/4} x \left(x^2+1\right)^{1/4}-\xi_1  \left(\xi_1^2+1\right)^{1/4} \left(x^2+1\right)^{3/4}-(x-\xi_1),
\end{equation*}
thus $f(\xi_1)=0$. Computing the derivative we have
\begin{equation*}
    f'(x)=\left(\xi_1 ^2+1\right)^{3/4} \left(x^2+1\right)^{1/4}+\frac{\left(\xi_1 ^2+1\right)^{3/4} x^2}{2 \left(x^2+1\right)^{3/4}}-\frac{3 \xi_1  \left(\xi_1^2+1\right)^{1/4} x}{2 \left(x^2+1\right)^{1/4}}-1,
\end{equation*}
thus $f'(\xi_1)=0$. Computing the second-order derivative we have
\begin{equation*}
    f''(x)=\frac{3 \left(\xi_1^2+1\right)^{1/4} \left(x^2+2\right) \left(\sqrt{\xi_1 ^2+1} x-\xi_1  \sqrt{x^2+1}\right)}{4 \left(x^2+1\right)^{7/4}}.
\end{equation*}
Now let $g(x)=\frac{x}{\sqrt{1+x^2}}$, we see that $g'(x)=\frac{1}{(1+x^2)^{3/2}} \geq 0$ for all $x \in \mathbb{R}$. 
Hence
\begin{equation*}
    \frac{x}{\sqrt{1+x^2}} \geq \frac{\xi_1}{\sqrt{1+\xi_1^2}},\qquad x\geq\xi_1.
\end{equation*}
Thus for all $x \geq \xi_1$, we have $f''(x) \geq 0$, thereby $$f'(x) \geq f'(\xi_1) =0,$$
and furthermore,
$$f(\xi_2) \geq f(\xi_1)=0,\qquad \xi_2\geq\xi_1.$$
This proves Lemma \ref{lemA}.

\section{Proof of Proposition \ref{remark1}}
If we replace the exponent $3/4$ in \eqref{E2} by $\alpha$ so that\footnotemark\footnotetext{For the estimate \eqref{E2} to hold at $\xi_2=\infty$, we have to assume that $\alpha>1/2$.} $1/2<\alpha<3/4$
and allow the right hand side of \eqref{E2} multiplied by $1 \leq C < +\infty$, 
then to prove the resulted inequality we are led to consider, for fixed $\xi_1$, the following function 
 \begin{equation*}
     f(x)=C\left(x \left(\xi_1 ^2+1\right)^{\alpha} \left(x^2+1\right)^{\alpha-\frac{1}{2}}-\xi_1  \left(\xi_1 ^2+1\right)^{\alpha-\frac{1}{2}} \left(x^2+1\right)^{\alpha}\right)-(x-\xi_1).
 \end{equation*}
 Note that $f(\xi_1)=0$. Computing the derivative we have 
 \begin{align*}
     f'(x)=&C\Biggr(2 \alpha x \left(\xi_1 ^2+1\right)^{\alpha-\frac{1}{2}} \left(x^2+1\right)^{\alpha-\frac{3}{2}} \left(\sqrt{\xi_1 ^2+1} x-\xi_1  \sqrt{x^2+1}\right)\\
     &\qquad+\left(\xi_1 ^2+1\right)^{\alpha} \left(x^2+1\right)^{\alpha-\frac{3}{2}}\Biggr)
     -1.
 \end{align*}
Thus $$f'(\xi_1)=C\left(\xi_1 ^2+1\right)^{2 \alpha-\frac{3}{2}}-1.$$ Since $\alpha<3/4$ and $C < +\infty$, we have 
 $$\lim_{\xi_1\rightarrow +\infty} f'(\xi_1)=-1<0.$$ This means that if $\xi_1$ is sufficiently large and $\xi_2>\xi_1$ is near $\xi_1$, we would have 
 $$f(\xi_2) < f(\xi_1)=0,$$ namely 
 \begin{equation*}
     \frac{\xi_2-\xi_1}{(1+\xi_2^2)^{\alpha}(1+\xi_1^2)^{\alpha}} > 
    C\left(\frac{\xi_2}{\sqrt{1+\xi_2^2}} - \frac{\xi_1}{\sqrt{1+\xi_1^2}}\right), \quad 
    \xi_2 > \xi_1,
 \end{equation*}
 which contradicts the inequality for exponent $\alpha$. This proves Proposition \ref{remark1}.

\section{Proof of Theorem \ref{thm1}}

The proof of Theorem \ref{thm1} uses the elegant machinery of Ozawa and Rogers \cite{OzaRog14}.
Since we can not do better, for completeness we have to ``reload their arguments".

The Ozawa-Rogers machinery starts with writing
\begin{align*}
    e^{it\sqrt{1-\Delta}}&f_1(x)e^{it\sqrt{1-\Delta}}f_2(x) \\ 
    &=\frac{1}{(2\pi)^2} \int_{\mathbb{R}^{2}} \widehat{f_1}(\xi_1)\widehat{f_2}(\xi_2)
    e^{ix(\xi_1+\xi_2)+it\left(\sqrt{1+\xi_1^2}+\sqrt{1+\xi_2^2}\right)} d\xi_1 d\xi_2 \\
    &=\frac{1}{(2\pi)^2} \int_{\mathbb{R}^{2}} F(\xi_1,\xi_2) e^{ix(\xi_1+\xi_2)+it\left(\sqrt{1+\xi_1^2}+\sqrt{1+\xi_2^2}\right)} d\xi_1 d\xi_2,
\end{align*}
where $$F(\xi_1,\xi_2)=\frac{1}{2}\left(\widehat{f_1}(\xi_1)\widehat{f_2}(\xi_2)+\widehat{f_1}(\xi_2)\widehat{f_2}(\xi_1)\right).$$ By symmetry, we obtain
\begin{align*}
    e^{it\sqrt{1-\Delta}}&f_1(x)e^{it\sqrt{1-\Delta}}f_2(x) \\
    &= 
    \frac{2}{(2\pi)^2} \int_{\xi_2 \geq \xi_1} F(\xi_1,\xi_2) e^{ix(\xi_1+\xi_2)+it\left(\sqrt{1+\xi_1^2}+\sqrt{1+\xi_2^2}\right)} d\xi_1 d\xi_2\\
    &= 
    \frac{2}{(2\pi)^2} \int_{\xi_2 \geq \xi_1} F(\xi_1,\xi_2) 
    e^{ix\eta_1+it\eta_2} \frac{1}{|J(\xi_1,\xi_2)|} d\eta_1 d\eta_2,
\end{align*}
where we used the following change of variables
\begin{equation*}
    \eta_1=\xi_1+\xi_2 \quad \text{and} \quad \eta_2=\sqrt{1+\xi_1^2}+\sqrt{1+\xi_2^2},
\end{equation*}
with the Jacobian $J$ given by
\begin{equation*}
    J(\xi_1,\xi_2)=\left | \begin{matrix}
\frac{\partial\eta_1}{\partial\xi_1}& \frac{\partial\eta_1}{\partial\xi_2} \\
\frac{\partial\eta_2}{\partial\xi_1}& \frac{\partial\eta_2}{\partial\xi_2} \\
\end{matrix} \right | = \frac{\xi_2}{\sqrt{1+\xi_2^2}} - \frac{\xi_1}{\sqrt{1+\xi_1^2}}, \quad\xi_2\geq\xi_1.
\end{equation*}
By Plancherel's theorem and reversing the change of variables, we have 
\begin{equation} \label{E4}
\begin{aligned}
   & \left\| e^{it\sqrt{1-\Delta}}f_1e^{it\sqrt{1-\Delta}}f_2 \right\|_{L_{t,x}^2(\mathbb{R}^{1+1})}^2 \\
    &\qquad= 
    \frac{4}{(2\pi)^2} \int_{\xi_2 \geq \xi_1} |F(\xi_1,\xi_2)|^2 \frac{1}{|J(\xi_1,\xi_2)|^2} d\eta_1 d\eta_2 \\
    &\qquad= \frac{4}{(2\pi)^2} \int_{\xi_2 \geq \xi_1} |F(\xi_1,\xi_2)|^2 \frac{1}{|J(\xi_1,\xi_2)|} d\xi_1 d\xi_2. 
\end{aligned}
\end{equation}

At this point we need the following elementary inequality. 

\begin{lem} \label{lemKGunit}
For $\xi_2 \geq \xi_1$,
\begin{equation}\label{E5}
    \frac{\xi_2}{\sqrt{1+\xi_2^2}} - \frac{\xi_1}{\sqrt{1+\xi_1^2}} \geq 
    1- \frac{(1,\xi_1)\cdot(1,\xi_2)}{\sqrt{1+\xi_1^2}\sqrt{1+\xi_2^2}}.
\end{equation}
\end{lem}

Proof of \eqref{E5} is deferred to next section.
Substituting \eqref{E5} into \eqref{E4}, we have
\begin{align*}
    &\left\| e^{it\sqrt{1-\Delta}}f_1e^{it\sqrt{1-\Delta}}f_2 \right\|_{L_{t,x}^2(\mathbb{R}^{1+1})}^2\\
    &\qquad\leq 
    \frac{4}{(2\pi)^2} \int_{\xi_2 \geq \xi_1} |F(\xi_1,\xi_2)|^2 
    \left(1-\frac{(1,\xi_1)\cdot(1,\xi_2)}{\sqrt{1+\xi_1^2}\sqrt{1+\xi_2^2}}\right)^{-1} d\xi_1 d\xi_2\\
    &\qquad=\frac{2}{(2\pi)^2} \int_{\mathbb{R}^2} |F(\xi_1,\xi_2)|^2 
    \left(1-\frac{(1,\xi_1)\cdot(1,\xi_2)}{\sqrt{1+\xi_1^2}\sqrt{1+\xi_2^2}}\right)^{-1} d\xi_1 d\xi_2.
\end{align*}
Since $\text{\rm supp}\,(\widehat{f_1}) \cap \text{\rm supp}\,(\widehat{f_2})=\emptyset$, we have
\begin{equation*}
    |F(\xi_1,\xi_2)|^2 = \frac{1}{4}\left(|\widehat{f_1}(\xi_1)|^2|\widehat{f_2}(\xi_2)|^2+|\widehat{f_1}(\xi_2)|^2|\widehat{f_2}(\xi_1)|^2\right),
\end{equation*}
therefore
\begin{align*}
    &\left\| e^{it\sqrt{1-\Delta}}f_1e^{it\sqrt{1-\Delta}}f_2 \right\|_{L_{t,x}^2(\mathbb{R}^{1+1})}^2\\
    &\qquad\leq \frac{1}{2(2\pi)^2} \int_{\mathbb{R}^2} \left(|\widehat{f_1}(\xi_1)|^2|\widehat{f_2}(\xi_2)|^2+|\widehat{f_1}(\xi_2)|^2|\widehat{f_2}(\xi_1)|^2\right) \\
    &\qquad\qquad\qquad\times\left(1-\frac{(1,\xi_1)\cdot(1,\xi_2)}{\sqrt{1+\xi_1^2}\sqrt{1+\xi_2^2}}\right)^{-1} d\xi_1 d\xi_2\\
    &\qquad= \frac{1}{(2\pi)^2} \int_{\mathbb{R}^{2}} |\widehat{f_1}(\xi_1)|^2|\widehat{f_2}(\xi_2)|^2 \left(1-\frac{(1,\xi_1)\cdot(1,\xi_2)}{\sqrt{1+\xi_1^2}\sqrt{1+\xi_2^2}}\right)^{-1} d\xi_1 d\xi_2.
\end{align*}
Theorem \ref{thm1} is proved.

\section{Proof of Lemma \ref{lemKGunit}}

For fixed $\xi_1$ we consider 
\begin{equation*}
    f(x)=-\sqrt{\xi_1 ^2+1} \sqrt{x^2+1}-\xi_1  \sqrt{x^2+1}+\sqrt{\xi_1 ^2+1} x+\xi_1  x+1,
\end{equation*}
thus $f(\xi_1)=0$. Computing the derivative we have 
\begin{equation*}
    f'(x)=\frac{(\sqrt{\xi_1 ^2+1}+\xi_1 ) \left(\sqrt{x^2+1}-x\right)}{\sqrt{x^2+1}}.
\end{equation*}
Note that $f'(\xi_1)=\frac{1}{\sqrt{\xi_1 ^2+1}}$. Computing the second-order derivative we have 
\begin{equation*}
    f''(x)=\frac{-\sqrt{\xi_1 ^2+1}-\xi_1 }{\left(x^2+1\right)^{3/2}} \leq 0
\end{equation*}
for all $x \in \mathbb R$. Since
$\lim_{x\rightarrow +\infty} f'(x)=0$, we have $f'(x) \geq 0$ for all $x \geq \xi_1$. Hence $$f(\xi_2) \geq f(\xi_1)=0,\qquad\xi_2 \geq \xi_1.$$
This proves Lemma \ref{lemKGunit}.

\section{On the Jacobian estimates of Ozawa and Rogers} \label{sec:restr}

For $\xi_1,\xi_2 \in\bR^*=\bR\cup\{\pm\infty\}$, define
$$\chi(\xi_1,\xi_2)=\frac{|\xi_1-\xi_2|}{\sqrt{1+\xi_1^2}\sqrt{1+\xi_2^2}}$$
and
$$J(\xi_1,\xi_2)= \left|\frac{\xi_1}{\sqrt{1+\xi_1^2}}-\frac{\xi_2}{\sqrt{1+\xi_2^2}}\right|.$$
Here $\chi$ is the chordal distance (on the extended complex plane $\wh\bC$), 
and it is used for example in defining the spherical derivative of meromorphic functions (see Zalcman \cite{Zal98}).
As before, $J$ is the Jacobian determinant for the change of variables
$$ \eta_1=\xi_1+\xi_2 \quad \text{and} \quad \eta_2=\sqrt{1+\xi_1^2}+\sqrt{1+\xi_2^2}.$$
We have the following straightforward upper bounds: $\forall\,\xi_1,\xi_2 \in\bR^*$,
$$\chi(\xi_1,\xi_2)\leq1\quad\text{and}
\quad J(\xi_1,\xi_2)\leq2.$$

Following Ozawa and Rogers \cite{OzaRog14} we introduce

\begin{defn}
For $\xi_1,\xi_2 \in\bR^*$ and $\alpha\in[1,2]$, let
$$\sigma_\alpha(\xi_1,\xi_2)=\frac{|\xi_1-\xi_2|^\alpha}{(1+\xi_1^2)^{\frac12+\frac{\alpha}{4}}(1+\xi_2^2)^{\frac12+\frac{\alpha}{4}}}.$$
\end{defn}

Note that $\sigma_\alpha(\xi_1,\xi_2)=0$ when $\xi_1$ agrees with $\xi_2$.
For our purpose below, 
we also impose the constraint $\alpha\geq1$ on $\sigma_\alpha$ (namely, continuity of order $\geq1$).
As $\alpha\leq2$,
$$\begin{aligned}
\sigma_\alpha(\xi_1,\xi_2)&=[\chi(\xi_1,\xi_2)]^\alpha\frac{1}{(1+\xi_1^2)^{\frac12-\frac{\alpha}{4}}(1+\xi_2^2)^{\frac12-\frac{\alpha}{4}}}\\
&\leq[\chi(\xi_1,\xi_2)]^\alpha\leq1,
\end{aligned}$$
and in particular,
$$\sigma_2(\xi_1,\xi_2)=[\chi(\xi_1,\xi_2)]^2\leq1.$$
Thus $\{\sigma_\alpha\}_{\alpha\in[1,2]}$ can be interpreted as a scale of deformed chordal distances. 

The following natural extension of Lemma \ref{lemA}, which bound $J$ from below with $ \sigma_\alpha$,
was established by Ozawa and Rogers \cite{OzaRog14}, again via ``an unlikely combination of five trigonometric identities". 
These inequalities are crucial, and actually critical\footnotemark
\footnotetext{In their machinery, the power $\frac12+\frac{\alpha}{4}$ with $\alpha\in[1,2]$ gives the best possible restriction results.}, 
in proving their Fourier restriction results on the whole hyperbola.

\begin{lem}[Ozawa-Rogers] \label{lem:elem}
For all $\alpha\in[1,2]$ and all $\xi_1,\xi_2 \in\bR^*$,
\begin{equation} \label{eqn:elem}
  \sigma_\alpha(\xi_1,\xi_2)\leq 2^{\alpha-1}J(\xi_1,\xi_2).
\end{equation}
\end{lem}

\begin{rem}
The sharp constant $2^{\alpha-1}$ is implicit in \cite{OzaRog14}.
\end{rem}

\begin{rem}
J. Chang informed us a nice extension of \eqref{eqn:elem} for $\alpha=2$ to $\wh\bC$.
\end{rem}

The previous convexity proof for $\alpha=1$ seems to be inaffective at least when $\alpha=2$.
The aim of this section is to present a unfied approach to \eqref{eqn:elem} based on interpolation
and a simple yet very useful reformulation of the Jacobian determinant $J$.

\begin{proof}
We observe that the following interpolation (or, factorization) relation holds
$$\sigma_\alpha(\xi_1,\xi_2)=[\sigma_1(\xi_1,\xi_2)]^{2-\alpha}[\sigma_2(\xi_1,\xi_2)]^{\alpha-1}.$$
Thus it suffices to prove Lemma \ref{lem:elem} for $\alpha\in\{1, 2\}$.

Note that the case $\xi_1=\xi_2$ is trivial. The case $\xi_1=-\xi_2$ is also simple as
$$\begin{aligned}
\sigma_\alpha(-\xi_2,\xi_2)&=\frac{2^\alpha|\xi_2|^\alpha}{(1+\xi_2^2)^{1+\frac{\alpha}{2}}}\\
&\leq2^{\alpha-1}\frac{2|\xi_2|}{\sqrt{1+\xi_2^2}}=2^{\alpha-1}J(-\xi_2,\xi_2).
\end{aligned}$$
In above inequality, we used $1+\xi_2^2\geq1$, $\alpha\geq1$, and $|\xi_2|/\sqrt{1+\xi_2^2}\leq1$.

Now, let us first consider $\alpha=2$.
For $\xi_1\neq\pm\xi_2$, we start with reformulating $J$ as
\begin{equation} \label{eqn:Jredu}
\begin{aligned}
J(\xi_1,\xi_2)&=\left|\frac{\xi_1\sqrt{1+\xi_2^2}-\xi_2\sqrt{1+\xi_1^2}}{\sqrt{(1+\xi_1^2)(1+\xi_2^2)}} \right|\\
&=\frac{|\xi_1^2-\xi_2^2|}{\left|\xi_1\sqrt{1+\xi_2^2}+\xi_2\sqrt{1+\xi_1^2} \right|\sqrt{(1+\xi_1^2)(1+\xi_2^2)}}.
\end{aligned}
\end{equation}
Then, we have
$$\begin{aligned}
\sigma_2(\xi_1,&\xi_2)[J(\xi_1,\xi_2)]^{-1}\\
&=\frac{|\xi_1-\xi_2|}{|\xi_1+\xi_2|}
\frac{\left|\xi_1\sqrt{1+\xi_2^2}+\xi_2\sqrt{1+\xi_1^2}\right|}{\sqrt{(1+\xi_1^2)(1+\xi_2^2)}}.
\end{aligned}$$

\textbf{Case 1}. If $\xi_1$ and $\xi_2$ have same sign, then
\begin{equation} \label{eqn:samesign}
\begin{aligned}
\sigma_2(\xi_1,&\xi_2)[J(\xi_1,\xi_2)]^{-1}\\
&=\frac{|\xi_1-\xi_2|}{|\xi_1+\xi_2|}
\left|\frac{\xi_1}{\sqrt{1+\xi_1^2}}+\frac{\xi_2}{\sqrt{1+\xi_2^2}}\right|\leq2.
\end{aligned}
\end{equation}

\textbf{Case 2}. If $\xi_1$ and $\xi_2$ have opposite sign, then
$$\begin{aligned}
\sigma_2(\xi_1,&\xi_2)[J(\xi_1,\xi_2)]^{-1}\\
&=\frac{\left|\xi_1\sqrt{1+\xi_2^2}+\xi_2\sqrt{1+\xi_1^2}\right|}{|\xi_1+\xi_2|}
\frac{|\xi_1-\xi_2|}{\sqrt{(1+\xi_1^2)(1+\xi_2^2)}}\\
&=\frac{|\xi_1-\xi_2|}{\left|\xi_1\sqrt{1+\xi_2^2}-\xi_2\sqrt{1+\xi_1^2}\right|}
\frac{|\xi_1-\xi_2|}{\sqrt{(1+\xi_1^2)(1+\xi_2^2)}}\leq1.
\end{aligned}$$
In the last step we used $\chi(\xi_1,\xi_2)\leq1$ and the following elementary inequality
$$\frac{|A+B|}{\left|A\sqrt{1+B^2}+B\sqrt{1+A^2}\right|}\leq1,\quad A\geq0,\, B\geq0.$$

Next, we consider $\alpha=1$.
Using the reformulation \eqref{eqn:Jredu}, it suffices to prove
\begin{equation} \label{eqn:alpha1redu}
\left|\xi_1\frac{(1+\xi_2^2)^{\frac14}}{(1+\xi_1^2)^{\frac14}}+\xi_2\frac{(1+\xi_1^2)^{\frac14}}{(1+\xi_2^2)^{\frac14}}\right|\leq|\xi_1+\xi_2|.
\end{equation}
Without loss of generality we can assume $|\xi_2|\geq|\xi_1|$. By setting
$$s=\xi_1,\quad t=\xi_2, \quad \beta=\frac{(1+\xi_2^2)^{\frac14}}{(1+\xi_1^2)^{\frac14}},$$
\eqref{eqn:alpha1redu} becomes
$$s^2\beta^2+\frac{t^2}{\beta^2}\leq s^2+t^2,$$
which is further equivalent to
$$(s^2\beta^2-t^2)(\beta^2-1)\leq0.$$
The last inequality is valid since
$$1\leq\beta^2\leq\frac{|\xi_2|}{|\xi_1|}\leq \frac{\xi_2^2}{\xi_1^2}=\frac{t^2}{s^2},$$
the second estimate involved being the monotonicity of $\xi/\sqrt{1+\xi^2}$.

This proves Lemma \ref{lem:elem} (the constant $2^{\alpha-1}$ also arrives under interpolation).
\end{proof}

\begin{rem}
In \eqref{eqn:elem} for $\alpha=1$, the constant 1 is sharp as we can take $\xi_2\rightarrow\xi_1=0$.
For $\alpha=2$ in \eqref{eqn:elem}, the constant 2 is also sharp as we can take $\xi_2>>\xi_1\rightarrow\infty$ in \eqref{eqn:samesign}.
\end{rem}

\begin{rem}
Using the reformulation \eqref{eqn:Jredu} again, in order that 
$$\sigma_\alpha(\xi_1,\xi_2)[J(\xi_1,\xi_2)]^{-1}\leq C_\alpha\quad\text{as}\quad\xi_2\rightarrow\xi_1,$$
with some bound $C_\alpha<\infty$, it is necessary that $\alpha\geq1$.
\end{rem}

\section*{Acknowledgement}

Research of the authors is supported by the National Natural Science Foundation of China (Grant No. 11801274)
and the Jiangsu Province Natural Science Foundation (Grant No. BK20200308).
The second author (YCH) thanks Professors Jianming Chang (CSLG) and Yan Xu (NJNU) for helpful communications on chordal distance.
The authors thank sincerely an anonymous referee for a kind suggestion.

\bigskip

\section*{\textbf{Compliance with ethical standards}}

\bigskip

\textbf{Conflict of interest} The authors have no known competing financial interests or
personal relationships that could have appeared to influence this reported work.

\bigskip

\textbf{Availability of data and material} Not applicable.

\bigskip

\bibliographystyle{alpha}
%\bibliography{main}

\end{document}